\documentclass [12pt,a4paper,reqno]{amsart} \textwidth 165mm
\usepackage{amsmath,amsthm}
\usepackage{amsfonts,amssymb}
\allowdisplaybreaks

\newcommand{\be}{\begin{equation}}
\newcommand{\ef}{\end{equation}}
\chardef\bslash=`\\ % p. 424, TeXbook
\newtheorem{thm}{Theorem}[section]
\newtheorem*{thm*}{Theorem}

\newtheorem{lem}[thm]{Lemma}
\newtheorem{prop}[thm]{Proposition}

\theoremstyle{definition}

\newtheorem*{remark*}{Remarks}
\newtheorem*{defn*}{Definition}

\theoremstyle{remark}

\numberwithin{equation}{section}

\newcommand{\G}{\Gamma}
\newcommand{\wt}{\widetilde}
\newcommand{\wh}{\widehat}
\newcommand{\fc}{\frac}
\newcommand{\bk}{\bigskip}
\newcommand{\iy}{\infty}
 \renewcommand{\sectionmark}[1]{}

\renewcommand{\Re}{\operatorname{Re}}
\newcommand{\Be}{Beltrami}
\newcommand{\hol} {holomorphic}
\newcommand{\qc} {quasiconformal}
\newcommand{\sh} {subharmonic}

\newcommand{\ve}{\varepsilon}
\newcommand{\e}{\epsilon}

\newcommand{\Te} {Teichm\"{u}ller}
\newcommand{\Ko} {Kobayashi}
\newcommand{\Ca} {Carath\'{e}odory}
\newcommand{\uTs} {universal Teichm\"{u}ller space}

\newcommand{\const}{\operatorname{const}}

 \usepackage{amsfonts}
\newcommand{\field}[1]{\mathbb{#1}}
\newcommand{\g}{\gamma}

\newcommand{\dl}{\delta}
\newcommand{\D}{\Delta}
\newcommand{\om}{\omega}
\newcommand{\z}{\zeta}
\newcommand{\ov}{\overline}
\newcommand{\vp}{\varphi}
\newcommand{\hC}{\wh{\field{C}}}
\newcommand{\C}{\field{C}}
\newcommand{\R}{\field{R}}

\newcommand{\B}{\mathbf{B}}
\newcommand{\T}{\mathbf{T}}

\newcommand{\Hol}{\operatorname{Hol}}

\newcommand{\x} {\mathbf x}

\renewcommand{\a} {\alpha}
\newcommand{\vr} {\varrho}
\newcommand{\ld}{\lambda}
\newcommand{\kp}{\kappa}

\newcommand{\Pot}{\operatorname{Pot}}

\newcommand{\Mob}{\operatorname{Mob}}

\begin{document}

\title{Hyperbolic geodesics, Krzyz's conjecture and beyond}
\author{Samuel L. Krushkal}

\begin{abstract} In 1968, Krzyz conjectured that for non-vanishing \hol \
functions  $f(z) = c_0 + c_1 z + \dots$ in the unit disk with
$|f(z)| \le 1$, we have the sharp bound $|c_n| \le 2/e$ for all $n
\ge 1$, with equality only for the function $f(z) = \exp [(z^n -
1)/(z^n + 1)]$ and its rotations. This conjecture was considered by
many researchers, but only partial results have been established.
The desired estimate has been proved only for $n \le 5$.

We provide here two different proofs of this conjecture and its
generalizations based on completely different ideas.
\end{abstract}

%\date{\today\hskip4mm({twoproofs.tex})}

\maketitle

\bigskip

{\small {\textbf {2000 Mathematics Subject Classification:} Primary:
30C50, 30C55, 32Q45; Secondary 30F45}

\medskip


\bigskip

\markboth{Samuel L. Krushkal}{Hyperbolic geodesics, Krzyz's
conjecture and beyond} \pagestyle{headings}

\bk

\section{Introduction, statement of results}.

\subsection{Conjecture}
Non-vanishing \hol \ functions $f(z) = c_0 + c_1 z + ...$ on the
unit disk $\D = \{z : |z| < 1\}$ (i.e., such that $f(z) \ne 0$ in
$\D$) form the normal families admitting certain invariance
properties, for example, the invariance under action of the
M\"{o}bius group of conformal self-maps of $\D$, complex
homogeneity, etc. One of the most interesting examples of such
families is the set $\mathcal B_1 \subset H^\iy$ of \hol \ maps of
$\D$ into the punctured disk $\D_{*} = \D \setminus \{0\}$.

Compactness of $\mathcal B_1$ in topology of locally uniform
convergence on $\D$ implies the existence for each $n \ge 1$ the
extremal functions $f_0$ maximizing $|c_n(f)|$ on $\mathcal B_1$.
Such functions are nonconstant and must satisfy $|f(e^{i \theta})| =
1$ for almost all $\theta \in [0, 2 \pi]$.

Estimating coefficients on $\mathcal B_1$ was originated in 1940's
(see \cite{R}). In 1968, Krzyz \cite{Kz} conjectured that for all $n
\ge 1$,
 \be\label{1.1}
|c_n| \le 2/e,
\end{equation}
with equality only for the function $\kp_n(z) = \kp(z^n)$, where
 \be\label{1.2}
\kp(z) := \exp \Bigl(\frac{z -1}{z + 1} \Bigr) = \fc{1}{e} +
\fc{2}{e} z - \fc{2}{3e} z^3 + ... \ .
\end{equation}
and its rotations $\epsilon_1 \kp(\epsilon_2 z)$ with $|\epsilon_1|
= |\epsilon_2| = 1$. Note that $\kp(z)$ is a \hol \ universal
covering map $\D \to \D_{*}$ moving $0$ to $1/e$.

This problem has been investigated by a large number of
mathematicians, however it still remained open. The estimate (1.1)
was established only for some initial coefficients $c_n$ including
all $n \le 5$ (see \cite{HSZ}, \cite{PS}, \cite{Sa}, \cite{Sz},
\cite{Ta}). On developments related to this problem see, e.g.,
\cite{Ba}, \cite{Ho}, \cite{HSZ}, \cite{LS}, \cite{MSUV}, \cite{R},
 \cite{Sz}.

Our main goal is to prove that Krzyz's conjecture is true for all $n
\ge 1$:

\begin{thm} For every $f(z) = c_0 + c_1 z + ... \in \mathcal B_1$
and $n \ge 1$, we have the sharp bound (1.1), and the equality
occurs only for the function $\kp_n$ and its rotations.
\end{thm}

\subsection{Proofs and generalizations}
We provide two completely different proofs of this theorem. The
first proof relies on complex geometry of convex Banach domains and
reveals geodesic features of the cover function (1.2), while the
second one involves the results related to the universal
Teichm\"{u}ller space and extremal Beltrami coefficients following
the lines originated in \cite{Kr2}. Both proofs shed light on the
intrinsic connection between the complex geodesics and extremals of
holomorphic functionals.

We also obtain some generalizations of Theorem 1.1. The arguments in
the first proof of Theorem 1.1 yields in the case $n = 2$ the
following improvement of this theorem: {\em any functional $J(f) =
c_2 + P(c_1)$, where $P(c_1)$ is a homogeneous polynomial of degree
$2$, satisfying $|P(c_1)| < \kp^\prime(0) = 2/e$ for all $f \in
\mathcal B_1$ distinct from $\kp$, satisfies
 \be\label{1.3}
\max_{\mathcal B_1} |J(f)| = \max_{\mathcal B_1} |c_2| = 2/e
\end{equation}
with the same extremal function $\kp_2$ (up to rotations)}.

The second proof deals with more general bounded functionals $J(f) =
c_n + F(c_{m_1}, \dots, c_{m_s}) $ on $\mathcal B_1$ where $c_j =
c_j(f); 1 \le n, m_j$ and $F$ is a holomorphic function of $s$
variables in an appropriate domain of $\C^s$. Assuming that this
domain contains the origin $\mathbf 0$ and that $F, \ \partial F$
vanish at $\mathbf 0$, we establish that any such a functional is
sharply estimated on appropriate subsets $E_r \subset \mathcal B_1$
(with $r$ depending on $n$) by
$$
\max_{E_r} |J(f)| = \max_{E_r} |c_n| = M_n r, \quad M_n =
\max_{\mathcal B_1} |J(f)|
$$
and obtain the desired bound (1.1) in the limit as $r \to 1$.

\bk
\section{Background: Invariant metrics on convex Banach domains}

We present briefly some basic results in complex geometry of convex
domains in complex Banach spaces, underlying the proofs of our main
theorem.

\subsection{Invariant metrics}
Let $M$ be a complex Banach manifold modeled by a Banach space $X$.
The  {\bf \Ko \ metric} $d_M$ on $M$ is the largest pseudometric $d$
on $M$ that does not get increased by \hol \ maps $h: \ \D \to M$ so
that for any two points $\x_1, \ \x_2 \in M$, we have
$$
d_M(\x_1, \x_2) \le \inf \{d_\D(0,t): \ h(0) = \x_1, \ h(t) =
\x_2\},
$$
where $d_\D$ is the hyperbolic metric on the unit disk of Gaussian
curvature $- 4$, hence with the differential form
$$
ds = \ld_\D(z)|dz| : = |dz|/(1 - |z|^2).
$$
The {\bf \Ca} \ distance between $\x_1$ and $\x_2$ in $M$ is
$$
c_M(\x_1, \x_2) = \sup d_\D(f(\x_1), f(\x_2)),
$$
where the supremum is taken over all \hol \ maps $f: \ E \to \D$.

In the case of a bounded domain $M$, both $d_M$ and $c_M$ are
distances (i.e., separate the points in $M$). For general properties
of invariant metrics we refer to \cite{Di}, \cite{Ko}. A remarkable
fact is:

\begin{prop} \cite{DTV}, \cite{Le} If $M$ is a convex domain in complex
Banach space, then
 \be\label{2.1}
d_M(\x_1, \x_2) = c_M(\x_1, \x_2) = \inf \{d_\D (h^{-1}(\x_1),
h^{-1}(\x_2)): \ h \in \Hol(\D, M)\}.
\end{equation}
\end{prop}

Similar equality holds for the differential (infinitesimal) forms of
these metrics which are defined on the tangent bundle $\mathcal T M$
of $M$.

\subsection{Complex geodesics}
A \hol \ map $h$ of the disk $\D$ into a Banach manifold endowed
with a pseudo-distance $\rho$ is called {\bf complex
$\rho$-geodesic} if there exist two points $t_1 \ne t_2 \in \D$ for
which
$$
d_\D(t_1, t_2) = \rho(h(t_1), h(t_2))
$$
(one says also that their images $h(t_1)$ and $h(t_2)$ can be joined
in $M$ by a complex $\rho$-geodesic; cf. \cite{Ve}). Any
$c_M$-geodesic is also $d_M$-geodesic, and then the equality (2.1)
holds for all points of the disk $h(\D)$.

Certain conditions ensuring the existence of complex geodesics,
which will be used here, are given in \cite{Di}, \cite{DTV}.

Assume that a Banach space $X$ has a predual space $Y$, i.e., that
$X = Y^\prime$ is the space of bounded linear functionals $x(y) =
<x, y>$ on $Y$, and consider on $X$ the {\bf weak$^*$ topology}
$\sigma(X, Y)$ which is the topology of pointwise convergence on
points of $Y$, i.e., a sequence $\{x_n\} \subset X$ is convergent in
$\sigma(X, Y)$ to $x \in X$ if $x_n(y) \to x(y)$ for all $y \in Y$.

If $X$ has a predual $Y$, then by the Alaoglu-Bourbaki theorem, the
closure of the open unit ball $X_1$ of the space $X$ in the topology
$\sigma(X, Y)$ is compact.

\bk

\begin{prop} \cite{Di}, \cite{DTV} Let $M$ be a bounded convex
domain in a complex Banach space $X$ with predual $Y$. If the
closure of $M$ is $\sigma(X, Y)$-compact, then every distinct pair
of points in $M$ can be joined by a complex $c_M$-geodesic.
\end{prop}

\bk This proposition also has its differential counterpart which
provides that under the same assumptions, for any point $\x \in M$
and any nonzero vector $v \in X$, there exists at least one complex
geodesic $h: \ \D \to M$ such that $h(0) = \x$ and $h^\prime(0)$ is
collinear to $v$ (cf. \cite{DTV}).

\subsection{Evaluation of \hol \ maps on geodesic disks}
We shall need the following corollary of the above propositions
controlling the growth of \hol \ maps with critical points on
geodesic disks (cf. \cite{Kr5}).

\begin{lem} Let a domain $M$ satisfy the assumptions of
Proposition 2.2 and $g$ be a \hol \ map $M \to \D$ whose restriction
to a geodesic disk $h(\D) \subset M, \ h(0) = \mathbf 0$, has at the
origin zero of order $m \ge 1$, i.e.,
$$
g \circ h(t) = c_m t^m + c_{m+1} t^{m+1} + \dots.
$$
Then the growth of $|g|$ on this disk is estimated by
 \be\label{2.2}
\begin{aligned}
|g \circ h(t)| &\le |t|^m (|t| + |c_m|)/(1 + |c_m||t|)  \\
&= \tanh d_M \Bigl(\mathbf 0, h \Bigl(|t|^m \ \fc{|t| + |c_m|}{1 +
|c_m||t|}\Bigr)\Bigr) \le \tanh d_M (\mathbf 0, h(t^m)).
\end{aligned}
\end{equation}
The equality in the right inequality occurs (even for one $t_0 \ne
0$) only when $|c_m| = 1$; then $h(t)$ is a hyperbolic isometry of
the unit disk and all terms in (2.2) are equal.
\end{lem}

\medskip\noindent
\textbf{Proof}. By Golusin's version of Schwarz's lemma, a \hol \
function
$$
f(t) = c_m t^m + c_{m+1} t^{m+1} + \dots: \D \to \D \quad (c_m \neq
0, \ \ m \ge 1)
$$
is estimated in $\D$ by
$$
|f(t)| \le |t|^m \fc{|t| + |c_m|}{1 + |c_m| |t|},
$$
and the equality occurs only for $f_0(t) = t^m (t+ c_m)/(1 + \ov c_m
t)$ (see [7, Ch. 8]).

It follows from Proposition 2.1 and weak$^*$ compactness of the
closure of $M$ in $\sigma(X, Y)$ that for any $t_0 \ne 0$ and $\x_0
= h(t_0)$ there exists a \hol \ map $j: \ M \to \D$ such that
$$
d_\D(0, j(\x_0)) = c_M (\mathbf 0, \x_0) = d_M (\mathbf 0, \x_0).
$$
Letting
$$
\eta(t) = |t|^m (|t|+ |c_m|)/(1 + |c_m| |t|),
$$
one gets $\eta(t) \le |t|$ and
$$
|g \circ h(t_0)| \le |j \circ h(\eta(t_0))| = \tanh d_M(\mathbf 0,
h(\eta(t_0)) \le \tanh  d_M(\mathbf 0, h(t_0)),
$$
which yields (2.2).

There is also a differential analog of the inequalities (2.2) which
involves the infinitesimal \Ca \ and \Ko \ metrics. It will not be
used here.

Lemma 2.3 straightforwardly extends to general complex Banach
manifolds $M$ having equal \Ca \ and \Ko \ distances.

\subsection{Generalized Gaussian curvature of subharmonic metrics}
The proof of Theorem 1.1 involves subharmonic conformal metrics
$\ld(t) |dt|$ on the disk having the curvature at most $- 4$ in a
somewhat generalized sense. As well-known, the Gaussian curvature of
a $C^2$-smooth metric $\ld > 0$ is defined by
$$
k_\ld (t) = - \fc{\D \log \ld(t)}{\ld(t)^2},
$$
where $D$ means the Laplacian $4 \partial^2/\partial z \partial \ov
z$.

A metric $\ld(t) |dt|$ in a domain $G \subset \C$ (or on a Riemann
surface) has the curvature less than or equal to $K$ {\bf in the
supporting sense} if for each $K^\prime > K$ and each $z_0$ with
$\ld(z_0) > 0$, there is a $C^2$-smooth supporting metric $\wh \ld$
for $\ld$ at $t_0$ (i.e., such that $\wt \ld(t_0) = \ld(t_0)$ and
$\wh \ld(t) \le \ld(t)$ in a neighborhood of $t_0$) with $k_{\wh
\ld}(t_0) \le K^\prime$, or equivalently, \be\label{2.3} \D \log \ld
\ge K \ld^2,
\end{equation}
A metric $\ld$ has curvature at most $K$ {\bf in the potential
sense} at $t_0$ if there is a disk $U$ about $z_0$ in which the
function
$$
\log \ld + K \Pot_U(\ld^2),
$$
where $\Pot_U$ denotes the logarithmic potential
$$
\Pot_U h = \fc{1}{2 \pi} \int\limits_U h(\z) \log |\z - t| d \xi d
\eta \quad (\z = \xi + i \eta),
$$
is \sh. This is equivalent to $\ld$ to satisfy (2.3) in the sense of
distributions.

One can replace above $U$ by any open subset $V \subset U$, because
the function $\Pot_U(\ld^2) - \Pot_V(\ld^2)$ is harmonic on $U$.

Due to Royden \cite{Ro2}, a conformal metric has curvature at most
$K$ in the supporting sense has curvature at most $K$ also in the
potential sense.

The following lemma concerns  the circularly symmetric (radial)
metrics on the disk (i.e. such that $\ld(t) = \ld(|t|)$) and is a
slight improvement of the corresponding Royden's lemma \cite{Ro2} to
singular metrics with a prescribed singularity at the origin.

\begin{lem} \cite{Kr4} Let $\ld(|t|) d |t|$ be a circularly symmetric
\sh \ metric on $\D$ such that
 \be\label{2.4}
\ld(r) = m c r^{m-1} + O(r^{m}) \quad \text{as} \ \ r \to 0 \ \
\text{with} \ \ 0 < c \le 1 \ \ (m = 1, 2, \dots),
\end{equation}
and this metric has curvature at most $- 4$ in the potential sense.
Then
 \be\label{2.5}
\ld(r) \ge \fc{m c r^{m-1}}{1 - c^2 r^{2m}}.
\end{equation}
\end{lem}

Note that all metrics subject to (2.4) are dominated by $\ld_m(t) =
m |t|^{m-1}/(1 - |t|^{2m})$.

\bk

\section{Preliminary results}

We first establish some analytic and geometric facts for
nonvanishing functions essentially applied in the proofs. These
results have their intrinsic interest.

\medskip
\noindent {\bf 1. Covering maps}.

\begin{prop} {\em (a)} \ Every function $f \in \mathcal B_1$
admits factorization \be\label{3.1} f(z) = \kp \circ \wh f(z),
\end{equation}
where $\wh f$ is a \hol \ map of the disk $\D$ into itself (hence,
from $H_1^\iy$) and $\kp$ is the function (1.2).

(b) \ Moreover, the map (3.1) generates an $H^\iy$-\hol \ map
$\mathbf k: \ H_1^\iy \to  \mathcal B_1$.
\end{prop}

\medskip
\noindent \textbf{Proof}. (a) Due to a general topological theorem,
any map $f: M \to N$, where $M, N$ are manifolds, can be lifted to a
covering manifold $\wh N$ of $N$, under appropriate relation between
the fundamental group $\pi_1(M)$ and a normal subgroup of $\pi_1(N)$
defining the covering $\wh N$ (see, e.g, [Ma]). This construction
produces a map $\wh f: M \to \wh N$ satisfying \be\label{3.2} f = p
\circ \wh f,
\end{equation}
where $p$ is a projection $\wh N \to N$. The map $\wh f$ is
determined up to composition with the covering transformations of
$\wh N$ over $N$ or equivalently, up to choosing a preimage of a
fixed point $x_0 \in \wh N$ in its fiber $p^{-1}(x_0)$. For \hol \
maps and manifolds the lifted map is also \hol.

In our special case, $\kp$ is a \hol \ universal covering map $\D
\to \D_{*} = \D \setminus \{0\}$, and the representation (3.2)
provides the equality (3.1) with the corresponding $\wh f$
determined up to covering transformations of the unit disk
compatible with the covering map $\kp$.

The assertion {\em (b)} is a consequence of a well-known property of
bounded \hol \ functions in Banach spaces with sup norm given by

\begin{lem} Let $E, \ T$ be open subsets of complex Banach spaces
$X, Y$ and $B(E)$ be a Banach space of \hol \ functions on $E$ with
sup norm. If $\vp(x, t)$ is a bounded map $E \times T \to B(E)$ such
that $t \mapsto \vp(x, t)$ is \hol \ for each $x \in E$, then the
map $\vp$ is \hol.
\end{lem}

Holomorphy of $\vp(x, t)$ in $t$ for fixed $x$ implies the existence
of complex directional derivatives
$$
\vp_t^\prime(x,t) = \lim\limits_{\z\to 0} \fc{\vp(x, t + \z v) -
\vp(x, t)}{\z} = \fc{1}{2 \pi i} \int\limits_{|\xi|=1} \fc{\vp(x, t
+ \xi v)}{\xi^2} d \xi,
$$
while the boundedness of $\vp$ in sup norm provides the uniform
estimate
$$
\|\vp(x, t + c \z v) - \vp(x, t) - \vp_t^\prime(x,t) c v\|_{B(E)}
\le M |c|^2,
$$
for sufficiently small $|c|$ and $\|v\|_Y$ (cf. \cite{Ha}).

The map $\mathbf k: \ \wh f \mapsto f$ is bounded on the ball
$H_1^\iy$. Applying Hartog's theorem on separate holomorphy to the
sums $g(z, t) = \wh f(z) + t \wh h(z)$ of $\wh f \in H_1^\iy, \ \wh
h \in H_1$ and $t$ from a region $B \subset \hC$ so that $g(z, t)
\in H_1^\iy$, one obtains that $g(z, t)$ are jointly \hol \ in both
variables $(z, t) \in \D \times B$. Thus the  restriction of the map
$\mathbf k_0$ onto intersection of the ball $H_1^\iy$ with any
complex line $L = \{\wh f + t \wh h\}$ is $H^\iy$-\hol, and hence
this map is \hol \ as the map $H_1^\iy \to  \mathcal B_1$, which
completes the proof of Lemma 3.1.

As an immediate corollary of this lemma, one gets the following
known estimate, which will be used here.

\begin{lem} For any $f \in \mathcal B_1$,
\be\label{3.3} |c_1| \le 2/e,
\end{equation}
with equality only for the rotations $e^{i \a_1} \kp(e^{i \a} z)$ of
$\kappa$ (in particular, these functions maximize $|c_1|$ among the
\hol \ covering maps $\D \to \D_{*}$).
 \end{lem}

\medskip
\noindent {\bf Proof}. Given $f \in \mathcal B_1$ distinct from
$\kp$, one may rotate its covering map $\wh f$ in (3.1) to get $\wh
f(0) = a$, where $0 < a < 1$. By Schwarz's  lemma, $|\wh
f^\prime(0)| \le 1 - |a|^2 < 1$; hence,
$$
|f^\prime(0)| = |\kp^\prime(a)| |\wh f^\prime(0)| < |\kp^\prime(a)|
= \fc{2 e^{(a-1)/(a+1)}}{(a + 1)^2} < \fc{2}{e},
$$
which implies (3.3).

\bk We shall also lift the functions $f \in \mathcal B_1$ to the
universal cover of $\D_{*}$ by the left half-plane $\C_{-} = \{w \in
\C: \ \Re w < 0\}$ using the map $\kp \circ \sigma^{-1} = \exp$,
where
$$
\sigma(z) = (z - 1)/(z + 1): \ \D \to \C_{-}.
$$
These lifts of  $f$ are reduced to choice of branches of $\log f(z)$
determined by the values of $\log f(0)$ in $\C_{-}$.

\bk \noindent {\bf 2. Open domain of nonvanishing functions and its
\hol \ embedding}.

Consider the annuli
$$
A_r = \{r < |z| < 1\}, \quad 0 < r < 1,
$$
exhausting the punctured disk $\D \setminus \{0\}$, and let
$\mathcal B_r$ be the subset of nonvanishing functions $f \in
\mathcal B_1$ sharing the values in $A_r$, that is,
$$
\mathcal B_r = \{f \in \mathcal B_1: \ f(\D) \subset A_r\};
$$
then $\mathcal B_r \subset \mathcal B_{r^\prime}$ if $r > r^\prime$.
Put
$$
\mathcal B_1^0 = \bigcup_r \mathcal B_r;
$$
this union is located in the unit ball $H_1^\iy$ of the space $H^\iy
= H^\iy(\D)$. It will be convenient to regard the free coefficients
$c_0(f)$ as the constant elements of $\mathcal B_1$.

The following lemma provides some needed topological properties of
these sets.

\begin{lem}
{\em (a)} \ For any $r \in (0, 1)$, every point of $\mathcal B_r$
has a neighborhood $U(f, \epsilon(r))$ in $H^\iy$, which contains
only the functions belonging to some $\mathcal B_{r_{*}}$, where $0
< r_{*} = r_{*}(r) \le r$.

{\em (b)} \ Each set $\mathcal B_r$ is path-wise connective in
$H_1^\iy$.
\end{lem}

It follows that {\em the union $\mathcal B_1^0$ is a domain in
$H_1^\iy$ filled by nonvanishing functions on $\D$.} In particular,
it contains all functions $f \in \mathcal B_1$ which are \hol \ and
nonvanishing on the closed disk $\ov D$.

\medskip
\noindent {\bf Proof}. To prove the assertion (a), assume the
contrary, i.e., that for some $r$ such $r_{*}$ does note exist. Then
there exist a function $f_0 \in \mathcal B_r$ and the sequences of
functions $f_n \in H_1^\iy$ convergent to $f_0$,
 \be\label{3.4}
\lim\limits_{n\to \iy} \|f_n - f_0\|_{H^\iy} = 0
\end{equation}
and of points $z_n \in \D$ convergent to $z_0, \ |z_0| \le 1$, such
that either $f_n (z_n) = 0 \ (n = 1, 2, \dots)$ or
$$
f_n (z_n) \ne 0, \quad \text{but} \ \ \lim\limits_{n\to \iy} f_n
(z_n) = 0.
$$
The first case means that $f_n$ vanish in $\D$; in the second one,
we have a sequence of nonvanishing functions $f_n$ belonging to
different sets $\mathcal B_{r_n}$, which are indexed by $r_n \to 0$.

If $|z_0| < 1$, we immediately reach a contradiction, because then
the uniform convergence of $f_n$ on compact sets in $\D$ implies
$f_0(z_0) = 0$, which is impossible.

Let $|z_0| = 1$. The values of $f_0$ on $\D$ must run in the annulus
$A_r$, thus $\inf_\D |f_0(z)| \ge r$. Hence, for $n \ge n_0$,
$$
|f_n(z_n) - f_0(z_n)| \ge \big\vert |f_0(z_n)| - |f_n(z_n)|\big\vert
\ge \fc{r}{2},
$$
and by continuity, there exists a neighborhood $\D(z_n, \dl_n) =
\{|z - z_n| < \dl_n\}$ of $z_n$ in $\D$, in which $|f_n(z) - f_0(z)|
> r/3$ for all $z$. This implies
$$
 \|f_n - f_0\|_{H^\iy} \ge
\sup_{\D(z_n, \dl_n)} |f_n(z) - f_0(z)| > \fc{r}{3}.
 $$
This inequality must hold for all $n \ge n_0$, contradicting (3.4).
The part (a) is proved.

To show that $\mathcal B_r$ is path-wise connective take its
arbitrary distinct points $f_1, \ f_2$. Similar to (3.1) one gets
$$
f_j = \chi_r \circ \wt f_j, \quad j = 1, 2,
$$
where $\wt f_j \in H_1^\iy$ and $\chi_r$ is a \hol \ universal
covering map $\D \to A_r$. Connecting the covers $\wt f_1$ and $ \wt
f_2$ in $H_1^\iy$ by the line interval $l_{1,2}(t) = t \wt f_1 + (1
- t)\wt f_2 \ 0 \le t \le 1$, one obtains a path $\chi_r \circ
l_{1,2}: \ [0, 1] \to \mathcal B_r$ connecting $f_1$ with $f_2$. The
continuity of $\chi_r \circ l_{1,2}$ in the norm of $H^\iy$ easily
follows from the fact that the covering map $\chi_r$ is reduced to
exponentiation (cf. Proposition 3.1). This completed the proof.

\bk Observe that this lemma does not contradict to existence of
sequences $\{f_n\} \in H_1^\iy$ of vanishing functions on $\D$ or of
$f_n$ with $\lim\limits_{n\to \iy} f_n(z_n) = 0$ convergent to $f_0
\in \mathcal B_1^0$ only uniformly on compact sets in $\D$.

Note also that $\mathcal B_1^0$ is dense in $\mathcal B_1$ in the
weak topology because any $f(z) = c_0 + c_1 z + \dots \in \mathcal
B_1$ is approximated locally uniformly, for example, by the homotopy
functions $f_r(z) = \om^{-1} [r \om(z)], \ 0 < r < 1$, with
$$
\om(z) = (z - c_0)/(1 - \ov c_0 z)
$$
mapping $\D$ holomorphically onto a subdomain $f_r(\D) \Subset \D$.
Hence,
$$
\sup_{\mathcal B_1^0} |c_n(f)| = \max_{\mathcal B_1} |c_n(f)|,
$$
and this supremum is attained only on $f_0 \in \mathcal B_1$ with
$\|f_0\|_\iy = 1$.

\bk Now take the branch of the logarithmic function $\log w = \log
|w| + i \arg w$ in the plane $\C_w$ slit along the positive real
semiaxes $\R_{+} = \{w = u + iv \in \C: \ u > 0\}$ for which $0 <
\arg w < 2 \pi$ (and hence $\log (- 1) = i \pi$).

Since $\D$ is simply connected and for every  $f \in \mathcal B_1^0$
we have $-\iy < \log |f(z)| < 0$ for all $z \in \D$, one can well
define the composition of $f$ with the chosen branch of the
logarithmic function, which generates a single valued \hol \
function
 \be\label{3.5}
\mathbf j_f(z) = \log f(z): \ \D \to \C_{-}.
\end{equation}
As was mentioned after Lemma 3.3, this means lifting $f$  to the
universal cover $\C_{-} \to \D \setminus \{0\}$ with the \hol \
universal covering map $\exp$.

Every such function $\mathbf j_f$ satisfies
 \be\label{3.6}
\sup_\D (1 - |z|^2)^\a |\log \mathbf j_f(z)| \le \sup_\D (1 -
|z|^2)^\a (\log |\mathbf j_f(z)| + |\arg \mathbf j_f(z)|) < \iy
\end{equation}
for any $\a > 0$. We embed the set $\mathbf j \mathcal B_1^0$ into
in the Banach space $\B$ of hyperbolically bounded \hol \ functions
on the disk $\D$ with norm
$$
\|\psi\|_\B = \sup_\D (1 - |z|^2)^2 |\psi(z)|.
$$
This space is dual to the space $A_1 = A_1(\D)$ of integrable \hol \
functions on $\D$ with $L_1$-norm, and every continuous linear
functional $l_\psi$ on $A_1$ can be represented by
  \be\label{3.7}
l_\psi(\vp) = \langle \psi, \vp\rangle_\D := \iint\limits_\D (1 -
|z|^2)^2 \ov{\psi(z)} \vp(z) dx dy
\end{equation}
with some $\psi \in \B$, uniquely determined by $l$ (see \cite{Be}).

We want to investigate the geometrical properties of the image
$\mathbf j \mathcal B_1^0$. First of all, we have

\begin{lem} The functions $\mathbf j_f \in \mathbf j \mathcal B_1^0$
fill a convex set in $\B$.
\end{lem}

\medskip
\noindent {\bf Proof}. Let $f_1, \ f_2$ be two distinct points in
$\mathcal B_1^0$; then their images $\psi_1 = \mathbf j f_1, \
\psi_2 = \mathbf j f_2$ are also different. The points of joining
interval $\psi_t = t \psi_1 + (1 - t) \psi_2$ with $0 \le t \le 1$
represent the functions $\mathbf j f_t = \log (f_1^t f_2^{1-t})$,
taking again the branch of logarithm defined above. For each $t$,
the product $f_1^t(z) f_2^{1-t}(z) \ne 0$ in $\D$, and $r <
|f_1(z)|^t |f_2(z)|^{1-t} < 1 - r$. Hence, this interval lies
entirely in $\mathbf j \mathcal B_1^0$.

\begin{lem} The map $\mathbf j$ is a \hol \ embedding of domain
 $\mathcal B_1^0$ into the space $\B$ carrying this domain
 onto a \hol \ Banach manifold modeled by $\B$.
\end{lem}

\medskip
\noindent {\bf Proof}. The map $\mathbf j: \ f \to \log f$ is
one-to-one, bounded on each subset $\mathcal B_r$ and continuous on
$\mathcal B_1^0$, which follows from Lemma 3.4 and (3.6).

To check its complex holomorphy, observe that each $f \in \mathcal
B_1^0$ belongs to subsets $\mathcal B_r$ with $r \le r_f$ (hence
$|f(z)| \ge r_f > 0$ in $\D$). Thus for any fixed $h \in H^\iy$ and
sufficiently small $|t|$ (letting $\mathbf j(f) = \mathbf j_f$),
$$
\mathbf j(f + th) - \mathbf j(f) = \log \Bigl( 1 + t \fc{h}{f}
\Bigr) = t \fc{h}{f} + O(t^2),
$$
with uniformly bounded remainder for $\|h\|_\iy \le c < \iy$. This
yields that the directional derivative of $\mathbf j$ at $f$ equals
$h/f$ and also belongs to $\B$.

In a similar way, one obtains that the inverse map $\mathbf j^{-1}:
\ \psi \to \e^{\psi}$ is \hol \ on intersections of a neighborhood
of $\psi$ in $\B$ with complex lines $\psi + t \om$ in $\mathbf j
\mathcal B_1^0$. The lemma is proved.

Both complex structures on $\mathbf j \mathcal B_1^0$ endowed by
norms on $H^\iy$ and on $\B$ are equivalent.

\bk \noindent {\bf 3. Complex geometry of sets $\mathbf j \mathcal
B_1^0$ and $\mathcal B_1^0$}.

\medskip
As a subdomain of a complex manifold modeled by $\B$, the set
$\mathbf j \mathcal B_1^0$ admits the invariant \Ko \ and \Ca \
metrics. Our goal is to show that the geometric features of this set
are similar to bounded convex domains in Banach spaces.

\begin{prop} (i) The \Ko \ and \Ca \ distances on
$\mathbf j \mathcal B_1^0$ are equal:
 \be\label{3.8}
d_{\mathbf j \mathcal B_1^0}(\psi_1, \psi_2) = c_{\mathbf j \mathcal
B_1^0} (\psi_1, \psi_2) = \inf \{d_\D (h^{-1}(\psi_1),
h^{-1}(\psi_2)): \ h \in \Hol(\D, \mathbf j \mathcal B_1)\},
\end{equation}
and similarly for the infinitesimal forms of these metrics.

(ii) Every two points in $\mathbf j \mathcal B_1^0$ can be joined by
$c$-geodesic (i.e., by a complex geodesic in the strongest sense).
\end{prop}

\medskip
\noindent {\bf Proof}. The equality (3.8) follows from the property
(ii). We establish this property in two steps.

$(a)$ \ First take the $\epsilon$-blowing up of $\mathbf j \mathcal
B_1^0$, that is, we consider the sets
$$
U_\epsilon = \bigcup_{\psi \in \mathbf j \mathcal B_1^0} \ \{\om \in
\B: \ \|\om - \psi \|_\B < \epsilon\}, \quad \epsilon > 0.
$$

For these sets, we have

\begin{lem} Every set $U_\epsilon$ is a (bounded) convex domain
in $\B$, and its weak$^*$-closure in $\sigma(\B, A_1)$ is compact.
\end{lem}

\medskip
\noindent {\bf Proof}. The openness and connectivity of $U_\epsilon$
are trivial. Let us check convexity. Take any two distinct points
$\om_1, \om_2$ in $U_\epsilon$ and consider the line interval
 \be\label{3.9}
\om_t = t \om_1 + (1 - t) \om_2, \quad 0 \le t \le 1,
\end{equation}
joining these points. Since, by definition of $U_\epsilon$, each
point  $\om_n \ (n = 1, 2)$ lies in the ball $B(\psi_n, \epsilon)$
centered at $\psi_n$ with radius $\epsilon$, and the interval
$\{\psi_t = t \psi_1 + (1 - t) \psi_2\}$ lies in $\mathbf j B_1^0$,
we have, for all $0 \le t \le 1$,
$$
\om_t - \psi_t = t(\om_1 - \psi_1) + (1 - t) (\om_2 - \psi_2)
$$
and
$$
\|\om_t - \psi_t\| \le t\|\om_1 - \psi_1\| + (1 - t) \|\om_2 -
\psi_2\| < \epsilon,
$$
which shows that the interval (3.9) lies entirely in $U_\epsilon$.

To establish $\sigma(\B, A_1)$-compactness of the closure $\ov
U_\epsilon$ in $\B$, note that weak$^*$ convergence of the functions
$\om_n \in \B$ to $\om$ implies the uniform convergence of these
functions on compact subsets of $\D$. It suffices to show that for
any bounded sequence $\{\om_n\} \subset \B$ we have the equality
 \be\label{3.10}
\lim\limits_{n\to \iy} \iint\limits_\D \fc{(1 - |\z|^2)^2
\om_n(\z)}{\z - z} d \xi d \eta = \iint\limits_\D \fc{(1 - |\z|^2)^2
\om(\z)}{\z - z} d \xi d \eta, \quad z \in \D^*,
\end{equation}
because the functions $w_z(\z) = 1/(\z - z)$ span a dense subset of
$A_1(\D)$. But if
$$
\sup_\D (1 - |\z|^2)^2 |\om(\z)| < M < \iy \quad \text{for all} \ \
n,
$$
the equality(3.10) is a consequence of Lebesgue's theorem on
dominant convergence. The lemma follows.

\bk $(b)$ \ We proceed to the proof of Proposition 3.7 and first
establish the existence of complex geodesics in domains $U_\epsilon,
\ \epsilon < \epsilon_0$. Our arguments follow \cite{DTV}.

Let $\om_1$ and $\om_2$ be distinct points in $U_\epsilon$. By
Proposition 2.1,
$$
d_{U_\epsilon}(\om_1, \om_2) = c_{U_\epsilon} (\om_1, \om_2) = \inf
\{d_\D (h^{-1}(\om_1), h^{-1}(\om_2)): \ h \in \Hol(\D,
U_\epsilon)\};
$$
hence there exists the sequences $\{h_n\} \subset \Hol(\D,
U_\epsilon)$  and $\{r_n\}, \ 0 < r_n < 1$, such that $h_n(0) =
\om_1$ and $h_n(r_n) = \om_2$ for all $n, \ \lim\limits_{n\to \iy}
r_n = r <1$ and $c_{U_\epsilon} (\om_1, \om_2) = d_\D(0, r)$. Let
$h_n(t) = \sum\limits_{m=0}^\iy a_{m,n} t^m$ for all $t \in \D$ and
$n$.

Take a ball $B(0, R) = \{\om \in \B: \ \|\om\| < R\}$ containing
$U_\epsilon$. For any $\om \in B(0, R)$, the Cauchy inequalities
imply $\|a_{n,m}\|_\B \le R$ for all $n$ and $m$. Passing, if
needed, to a subsequence of $\{h_n\}$, one can suppose that for a
fixed $m$, the sequence $a_{n,m}$ is weakly$^*$ convergent to $a_m
\in \B$ as $n \to \iy$, that is
$$
\lim\limits_{n\to \iy} \langle a_{n,m}, \vp\rangle_\D = \langle a_m,
\vp\rangle_\D \quad \text{for any} \ \ \vp \in A_1.
$$
Hence $h(t) = \sum\limits_{m=0}^\iy a_m t^m$ defines a \hol \
function from $\D$ into $\B$. Since $a_{n,0} = \om_1$ for all $n$,
we have $h(0) = \om_1$.

Now, let $\a, \ 0 <\a < 1$, and $\ve > 0$ be given. Choose $m_0$ so
that
$$
r \sum\limits_{m=m_0}^\iy \a^m < \ve.
$$
If $\vp \in A_1, \ \|\vp\| = 1$, then
$$
\sup_{|t|\le \a} | \langle h_n(t) - h(t), \vp \rangle_\D| \le
\sum\limits_{m=1}^{m_0-1} |\langle a_{n,m} - a_m, \vp \rangle_\D| +
2 r \sum\limits_{m=m_0}^\iy \a^m
$$
for all $n$, which implies that $h_n$ is convergent to $h$ in
$\sigma (\B, A_1)$ uniformly on compact subsets of $\D$ as $n \to
\iy$. Since $\ov U_\epsilon$ is $\sigma (\B, A_1)$ compact, $h(\D)
\subset \ov U_\epsilon$, and since $h(0)  \in U_\epsilon$, it
follows that $h(\D) \subset U_\epsilon$. For $r < r^\prime < 1$,
 \be\label{3.11}
\om_2 = h_n(r_n) = \fc{1}{2 \pi i} \int\limits_{|t| = r^\prime}
\fc{h_n(t)dt}{t - r_n} \to \fc{1}{2 \pi i} \int\limits_{|t| =
r^\prime} \fc{h(t)dt}{t - r} = h(r)
\end{equation}
as $n \to \iy$. Hence,
$$
d_\D(0, r) = c_{U_\epsilon} (\om_1, \om_2) = c_{U_\epsilon} (h(0),
h(r)),
$$
and $h$ is a $c$-geodesics in $U_\epsilon$.

There exists a \hol \ map $g: \D \to U_\epsilon$ such that for any
two points $t_1, t_2 \in \D$,
 \be\label{3.12}
d_\D(t_1, t_2) = d_{U_\epsilon}(g(t_1), g(t_2)) =
c_{U_\epsilon}(g(t_1), g(t_2)),
\end{equation}
and for any pair $(t, v), \ t \in \D, \ v \in \C$,
 \be\label{3.13}
\mathcal K_{U_\epsilon}(g(t), d g(t) v) = \fc{|v|}{1 - |t|^2}.
\end{equation}

\bk $(c)$ \ Let now $\om_1$ and $\om_2$ be two distinct points in
$\mathbf j \mathcal B_1^0$. Choose a decreasing sequence
$\{\epsilon_n\}$ approaching zero and take for every $n$ a complex
geodesic $h_n = h_{U_{\epsilon_n}}$ joining these points in
$U_{\epsilon_n}$, which was constructed in the previous step. Let
$g_n = g_{U_{\epsilon_n}}$ be the corresponding map $\D \to
U_{\epsilon_n}$ which provides the equalities (3.12), (3.13). Since
$d_\D$ is conformally invariant, one can take $g_n$ satisfying
$g_n^{-1}(\om_1) = 0, \ g_n^{-1}(\om_2) = r_n \in (0, 1)$. Then the
inequalities
$$
d_{U_{\epsilon_n}}(\om_1, \om_2) \le d_{U_{\epsilon_m}} (\om_1,
\om_2) \le d_{\mathbf j B_1} (\om_1, \om_2) \quad \text{for} \ \ m >
n
$$
imply $r_n \le r_m \le r_{*} < 1$, where $d_\D(0, r_{*}) =
d_{\mathbf j \mathcal B_1}(\om_1, \om_2)$. Hence, there exists
$\lim\limits_{n\to \iy} r_n = r^\prime \le r_{*}$.

The sequence $\{g_n\}$ is $\sigma(\B, A_1)$-compact and similar to
(3.11) the weak$^*$ limit of $g_n$ is a function $g \in \Hol(\D,
\mathbf j \mathcal B_1^0)$ which determines a complex geodesic for
both \Ko \ and \Ca \ distances on ${\mathbf j \mathcal B_1^0}$
joining the points $\om_1$ and $\om_2$ inside this set. Proposition
3.7 is proved.

An important consequence of Proposition 3.7 is that the initial
domain $\mathcal B_1^0$ in $H^\iy$ has similar complex geometric
properties, since the embedding $\mathbf j$ is biholomorphic. We
present it as

\begin{prop} (i) The \Ko \ and \Ca \ distances on domain
$\mathcal B_1^0$ and the corresponding infinitesimal metrics are
equal:
 \be\label{3.14}
\begin {aligned}
d_{\mathcal B_1^0}(f_1, f_2) &= c_{\mathcal B_1^0} (f_1, f_2) = \inf
\{d_\D (h^{-1}(f_1), h^{-1}(f_2)):
\ h \in \Hol(\D, \mathcal B_1)\}, \\
\mathcal K_{\mathcal B_1^0}(f, v) &= \mathcal C_{\mathcal B_1^0}(f,
v) \quad \text{for all} \ \ (f, v) \in T(\mathcal B_1^0).
\end{aligned}
\end{equation}

(ii) Every two points $f_1, f_2$ in $\mathcal B_1^0$ can be joined
by a complex geodesic.
\end{prop}

\bk
\section{First proof of Theorem 1.1}

This proof involves a complex homotopy of functions $f \in \mathcal
B_1^0$ and estimating the \Ko \ distance on the homotopy disks.

For any $f \in \mathcal B_1^0$, there is a complex \hol \ homotopy
connecting $f$ with $c_0(f)$ in $\mathcal B_1^0$. For $f = \kappa
\circ \wh f$ whose the cover $\wh f(z) = \wh c_1 z + \wh c_2 z^2 +
\dots \in H_1^\iy$ has zero free term, one can take
$$
\wh f_t(z) = \wh f(t z) = \wh c_1 t z + \dots: \ \D \times \D \to
\D_{*} \quad (\wh f_0(\cdot) = \mathbf 0)
$$
generating the underlying homotopy $f_t = \kp \circ \wh f_t$ in
$\mathcal B_1^0$.

In the case of generic $\wh f(z) = \wh c_0 + \wh c_m z + \dots \in
H_1^\iy \ (m \ge 1)$, we decompose it via
 \be\label{4.1}
\wh f = \om \circ \wh g_{\wh f},
\end{equation}
where
 \be\label{4.2}
 \begin{aligned}
g_{\wh f}(z) &= (\wh f(z) - \wh c_0)/(1 - \ov{\wh c_0} \wh f(z)), \\
\om(g) &= (\wh g + \wh c_0)/(1 + \ov{\wh c_0} \wh g),
\end{aligned}
\end{equation}
and set
 \be\label{4.3}
f_t(z) = \kp \circ \om(\wh g_{\wh f})(t z).
\end{equation}
The pointwise map $t \mapsto f_t$ generates by Lemma 3.2 a \hol \
map $\chi_f: \D \to H^\iy$, and the functional $J$ is
$n$-homogeneous with respect to this homotopy, $J(f_t) = t^n J(f)$.

Note also that for a fixed $\wh c_0$ (regarded again as a constant
function on $\D$), both maps in (4.2) are biholomorphic isometries
of the ball $H_1^\iy$; hence
 \be\label{4.4}
\om_{*} (\z) = \fc{\z \wh g_f/\|\wh g_f\|_\iy + \wh c_0}{1 + \ov{\wh
c_0}  \z \wh g_f/\|\wh g_f\|_\iy}: \ \D \to H_1^\iy
\end{equation}
carries out the complex geodesic $\z \mapsto \z \wh g_f/\|\wh
g_f\|_\iy$ into a complex geodesic in $H_1^\iy$ passing through $\wh
c_0$ and $\wh f$, and $\om_{*}(\z) = \wh f$ at $\z = \|\wh
g_f\|_\iy$.

By Proposition 3.9, there exists for each $f_t$ a complex geodesic
in $\mathcal B_1^0$ joining $f_t$ with $c_0(f)$; it determines a
\hol \ geodesic disk isometric to the hyperbolic plane $\mathbb
H^2$. We need to estimate the behavior of the distance $d_{\mathcal
B_1^0} ({f_t}, c_0)$ for $t \to 0$.

\begin{lem} Let $\wh f \in H_1^\iy$ have the expansion
$\wh f(z) = \wh c_m z^m + \dots$ with $\wh c_m \ne 0 \ ( m \ge 1)$,
and $\|\wh f\| = 1$. Then the geodesic parameter $\z$ and the
homotopy parameter $t$ are related near the origin via
 \be\label{4.5}
 |\z| = |\wh c_m| |t|^m + O(|t|^{m+1}), \quad t \to 0.
\end{equation}
\end{lem}

\medskip
\noindent \textbf{Proof}. Put $\wh p_m(z) = z^m$. The homotopy disk
$\D(\wh p_m)$ of this function in $H_1^\iy$ is filled by the
functions $\wh p_{m,t}(z) = t^m z^m$ with $|t| < 1$, while the
geodesic parameter on $\D(\wh p_m)$ is generated by hyperbolic
isometry $\z \mapsto \z \wh g/\|\wh g\|$). So, $\z = t^m$, and since
$$
\|\wh f_t - \wh c_m \wh p_{m,t}\|_{H^\iy} = |t|^{m+1} \|\wh c_{m+1}
+ t \wh c_{m+2} z + \dots \|_\iy = O(t^{m+1}),
$$
the relation (4.5) follows.

\bk
\begin{lem} For any $f = c_0 + c_1 z + \dots \in \mathcal B_1^0$,
we have the equality
 \be\label{4.6}
 d_{\mathcal B_1^0} (f, c_0) = \inf \{d_{H_1^\iy} (\wh f, \wh c_0):
 \ \kappa \circ \wh f = f\};
\end{equation}
moreover, there exists a map $\wh f^*(z) = c_0^* + c_1^* z + \dots$
covering $f$, on which the infimum in (4.6) is attained, i.e.,
 \be\label{4.7}
d_{\mathcal B_1^0} (f, c_0) = d_{H_1^\iy} (\wh f^*, \wh c_0^*).
\end{equation}
\end{lem}

\medskip
\noindent {\bf Proof}. We decompose the cover $\wh f(z) = \wh c_0 +
\wh c_1 z + \dots$ of $f$ in $H_1^\iy$ by (4.1), (4.2), getting
 \be\label{4.8}
d_{H_1^\iy}(\wh g_f, \mathbf 0) = d_\D (\|\wh g_f\|_\iy, 0) =
\tanh^{-1} (\|(\wh f - \wh c_0)/(1 - \ov{\wh c_0} \wh f)\|_\iy).
\end{equation}
and then apply to $\wh g_{\wh f}$ the transform (4.4). This yields a
complex geodesic in $H_1^\iy$ which connects $\wh c_0$ and $\wh f$.

Now observe that the universal covering map $\kappa_0: \D \to
\D_{*}$ extended by the equality (3.1) to all $\wh f \in H_1^\iy$
generates \hol \ map of the ball $H_1^\iy$ into domain $\mathcal
B_1^0$, which yields
$$
d_{\mathcal B_1^0}(f, c_0) = d_{\mathcal B_1^0} (\kappa_0 \circ \wh
f, \kappa_0(\wh c_0)) \le d_{H_1^\iy}(\wh f, \wh c_0),
$$
and
 \be\label{4.9}
d_{\mathcal B_1^0}(f, c_0) \le \inf_{\wh f} d_{H_1^\iy}(\wh f, \wh
c_0),
\end{equation}
where the infimum is taken over all covers $\wh f$ of $f$.

It remains to establish that in fact one has in (4.9) the equality
(so the infimum is attained). To prove this, assume to the contrary,
that
$$
d_{\mathcal B_1^0}(f, c_0) < \inf_{\wh f} d_{H_1^\iy}(\wh f, \wh
c_0),
$$
By Proposition 3.8, there exists a complex geodesic $h: \ \D \to
\mathcal B_1^0$ joining the points $c_0$ and $f$, and it follows
from the above,
$$
d_{\mathcal B_1^0}(f, c_0) = d_\D(\z_1, \z_2) < \tanh^{-1} \|(f -
c_0)/(1 - \ov c_0 f)\|_\iy,
$$
where $\z_1 = h^{-1}(c_0), \ \z_2 = h^{-1}(f)$. Lifting this $h$ by
(3.1) to its cover $\wh h$ of the unit disk into itself, one gets
the points $\wh h(\z_1), \wh h(\z_2)$ in $\D$ located in the fibers
over $c_0$ and $f$, respectively, and for these points
$$
d_{H_1^\iy}(\wh h(\z_1), \wh h(\z_2)) = d_\D(\z_1, \z_2) <
\tanh^{-1} \|(f - c_0)/(1 - \ov c_0 f)\|_\iy,
$$
which contradicts (4.8) and completes the proof of Lemma 4.2.

Consider the functional $J(f) = c_n$ on all $f \in H_1^\iy$. The map
$\mathbf k$ extending the function (1.2) holomorphically to
$H_1^\iy$ generates a \hol \ functional
$$
J_{\mathbf k} = J \circ \mathbf k: \ H_1^\iy \to \C,
$$
with $\max_{H_1^\iy} |J_{\mathbf k}| = \sup_{\mathcal B_1^0} |J|$.
We rescale this functional taking
$$
J_{\mathbf k}^0(\wh f) = \fc{J_{\mathbf k}(f)}{C_n}, \quad C_n =
\sup_{f\in \mathcal B_1^0} |J(f)| = \max_{H_1^\iy} |J_{\mathbf
k}(\wh f)|,
$$
which yields a \hol \ map of the ball $H_1^\iy$ onto the unit disk.
Similarly, let $J^0(f) = J(f)/C_n$.

We first estimate these functionals on the set of $f \in \mathcal
B_1^0$ whose covers $\wh f$ in $H_1^\iy$ are of the form
 \be\label{4.10}
\wh f(z) = \wh c_m z^m + \wh c_{m+1} z^{m+1} + \dots
\end{equation}
(with $\wh c_0 = 0, \ \wh c_m \ne 0 \ m \ge 1)$.

Fix a small $\vr > 0$ and let $|t| < \vr$. Take the geodesics $h_t:
\ \D \to H_1^\iy$ joining such $\wh f_t$ with $\mathbf 0$. Then
 \be\label{4.11}
J_{\mathbf k}^0(\wh f_t) = \wh \beta_n \z^n + \wh \beta_{n+1}
\z^{n+1} + \dots
\end{equation}
(where $\wh \beta_n \ne 0$ and all $|\wh \beta_j| \le 1$). Applying
Lemmas 2.3 and 4.1, one derives
$$
|J_{\mathbf k}^0(\wh f_t)| = |J_{\mathbf k}^0 \circ h_t(\cdot)| +
O(|t|^{mn+1}) \le |\wh \beta_n| |\wh c_m| |t|^{mn} + O(|t|^{mn+1}),
\quad t \to 0.
$$
This yields, similar to Lemma 4.2, the following upper bound for the
images $f_t = \kappa \circ \wh f_t \in \mathcal B_1^0$
 \be\label{4.12}
|J^0(f_t)| \le \inf |\wh \beta_n| |\wh c_m| |t|^{mn} + O(|t|^{mn+1})
= \inf |\wh \beta_n| |\wh c_m| |t|^{mn} + O(|t|^{mn+1})
\end{equation}
(each infimum again over $\wh f$ with $\kappa \circ \wh f = f$).

\bk Now we establish that the right-hand side of (4.12) yields
simultaneously the lower asymptotic bound for $|J^0(f_t)|$ with
small $|t|$, which means that (4.12) is in fact an asymptotic
equality.

\begin{lem} For any $f(z) = c_0 + c_m z^m + \dots \in \mathcal
B_1^0$ with $\wh f$ of the form (4.10), we have
 \be\label{4.13}
|J^0(f_t)| \ge \inf (|\wh \beta_n| |\wh c_m| |t|^{mn} +
O(|t|^{mn+1}),
\end{equation}
again taking the infimum over $\wh f$ with $\kappa \circ \wh f = f$.
\end{lem}

\bk\noindent{\bf Proof}. Fix again a small $\vr > 0$ and restrict
$J_{\mathbf k}^0$ to geodesic $h_t(\z) = h(\z; t)$ in $H_1^\iy$
joining $\mathbf 0$ with $\wh f_t, \ |t| \le \vr$. The corresponding
function $\wh g_t(\z) = J_{\mathbf k}^0 \circ h_t(\z)$ given by
(4.11) generates a conformal metric
$$
\ld_{\wh g}(\z; t) = \wh g_t^* \ld_\D(\z) = \fc{|\wh g^\prime(\z;
t)|}{1 - |\wh g(\z; t)|^2}
$$
of Gaussian curvature $- 4$ at noncritical points. The upper
envelope of these metrics
$$
\ld_{J^0}(\z; t) = \sup_{\wh f} \ld_{\wh g}(\z; t)
$$
(over the covers $\wh f$ of given $f$) is a \sh \ metric on $\D$
with curvature at most $- 4$ in the supporting sense and hence in
the potential sense. Averaging $\ld_{J^0}(\z; t)$ over the torus
$\{|\z| = |t| = r\}$ yields a circularly symmetric \sh \ metric
$$
\ld_{J^0}(r) = n |\wh \beta_n| r^n + O(r^{n+1})
$$
of curvature at most $-4$ in the potential sense. Estimating this
metric by Lemmas 2.4 and 4.1, one derives the following lower bound
 \be\label{4.14}
\ld_{J^0}(r) \ge \fc{m n \ c r^{mn-1}}{1 - c^2 r^{2mn}},
\end{equation}
where
$$
c = \inf_{\wh f} |\wh \beta_n||\wh c_m|
$$
(this estimate depends on $\vr$ which was fixed). On the other hand,
the hyperbolic length of the interval $[0, |\wh g_t(r)|]$ equals
$$
\tanh^{-1} |\wh g_t(r)| = \int\limits_0^{|\wh g_t(r)|} \fc{d \xi}{1
- |\xi|^2} = \int\limits_0^r \ld_{\wh g}(\xi; \cdot) d \xi,
$$
which implies (cf. \cite{Kr3})
  \be\label{4.15}
\tanh^{-1} |J^0(f_r)| = \sup_{\wh f} \int\limits_0^r \ld_{\wh
g}(\xi; \cdot) d \xi = \int\limits_0^r  \sup_{\wh f} \ld_{\wh
g}(\xi; \cdot) d \xi = \int\limits_0^r \ld_{J^0}(\xi) d \xi,
\end{equation}
where the second equality is obtained by taking a monotone
increasing subsequence of metrics
$$
\ld_1 = \ld_{\wh g_1}, \quad \ld_2 = \max (\ld_{\wh g_1}, \ld_{\wh
g_2}), \quad \ld_3 = \max (\ld_{\wh g_1}, \ld_{\wh g_2}, \ld_{\wh
g_3}), \dots
$$
corresponding to a sequence $\{\wh f_j\} \subset H_1^\iy$ for which
$\sup_j |J_{\mathbf k}^0(\wh f_j)| = \sup_{\mathcal B_1^0} |J^0(f)|$
and such that $\lim\limits_{j\to \iy} \ld_j = \sup_j \ld_{\wh g_j}$.
From (4.14) and (4.15),
$$
\int\limits_0^r \ld_{J^0}(\xi) d \xi \ge \tanh^{-1} (c r^{mn}) +
O(r^{n+1}), \quad r \to 0,
$$
which proves the desired estimate (4.13).

\bk We have established that for small $r > 0$,
 \be\label{4.16}
\max_{|t|=r} |J^0(f_t)| = r^{mn} |J^0(f)| = \max_{|t|=r} \inf_{\wh
f} |\wh \beta_n||c_m(f_t)| + O(r^{n+1}).
\end{equation}

\bk One can replace in the above arguments the cover $\kappa$ by any
universal covering map $\g^* \kappa = \kappa \circ \g: \ \D \to
\D_{*}$ with $\g \in \Mob (\D)$. Fix in (3.1) $\wh f$ with $\wh f(0)
= 0$ and choose $\g$ so that the point $\g(0) = \kappa^{-1}(c_0)$ is
placed in the closure of a fundamental triangle of a cyclic Fuchsian
group $\G$ representing $\D_{*}$ as the quotient $\D/\G$. Then,
instead of (3.1),
  \be\label{4.17}
f(z) = (\g^* \kp) \circ \wh f(z).
\end{equation}

Applying this to generic functions $f \in \mathcal B_1^0$ covered by
arbitrary $\wh f(z) = \wh c_0 + \wh c_m z^m + \dots \in H_1^\iy \ (m
\ge 1)$, one obtains similar to Lemma 4.2 that the equalities (4.6)
are extended to homotopy (4.3) as follows.

\begin{lem} For any
$\wh f(z) = \wh c_0 + \sum_m^\iy \wh c_n z^n \in H_1^\iy \ (\wh c_m
\ne 0, \ m \ge 1)$ and its image $f(z) = \kp \circ \wh f(z) = c_0 +
\sum_m^\iy c_n z^n \in \mathcal B_1^0 $,
 \be\label{4.18}
d_{\mathcal B_1^0} (f_t, c_0) = \inf_{\wh f} d_{H_1^\iy} (\wh f_t,
\mathbf 0) = \inf \{|\wh c_m(\wh f)|: \ \kp \circ \wh f = f\} |t|^m
+ O(|t|^{m+1}), \quad t \to 0,
\end{equation}
where each infimum is taken over the covers $\wh f$ of $f$ fixing
the origin and attained on some $\wh f$, and the estimate of the
remainder is in $H^\iy$-norm. Therefore,
 \be\label{4.19}
\max_{|t|=r} |J^0(f_t)| = \max_{|t|=r} \inf_{\wh f} |\wh
\beta_n||\wh c_m| r^{mn} + O(r^{mn+1}).
\end{equation}
\end{lem}

In particular, since $\wh \kp(z) = z$,
$$
d_{\mathcal B_1^0} (\kp_t, c_0) = |t| + O(|t|^2), \quad t \to 0.
$$
which shows that the {\em \hol \ disk $\D(\kappa)$ filled by the
homotopy functions $\kappa_t(z) = \kappa(tz), \ t \in \D$, is
geodesic in} $ \mathcal B_1^0$. Similarly, for $\kp_m(z) =
\kappa(z^m)$,
 \be\label{4.20}
d_{\mathcal B_1^0} (\kp_{m,t}, c_0) = |t|^m + O(|t|^{m+1}), \quad t
\to 0.
\end{equation}
In fact, the remainder terms in the last two equalities can be
omitted.

We can now complete the proof of Theorem 1.1. Let
$$
f^0(z) = c_0^0 + c_1^0 z + c_2^0 z^2 + \dots
$$
be an extremal function maximizing $|J(f)| \ (n > 1)$ on $\mathcal
B_1^0$; then $|J(f^0)| = C_n$. Rotating, if needed $f^0$, we get
$J^0(f^0) = 1$ and $J^0(f_r^0) = r$ for the homotopy  $f_r^0(z) =
f^0(r z), \ (0 < r <1)$.

Using this homotopy, we first show that $f^0$ must satisfy
 \be\label{4.21}
c_1^0 = 0.
\end{equation}
Indeed, assume that $c_1^0 \ne 0$ (hence, $C_n > J(\kp)$) and apply
Lemma 4.4 (with $m = 1$). By (4.18),
$$
d_{\mathcal B_1^0} (f_r^0, c_0^0) = |\wh c_1^0| r + O(r^2) =
\fc{|c_1^0| r}{|(\g^* \kp)^\prime(0)|} + O(r^2), \quad r \to 0,
$$
where $\wh c_1^0$ is the first coefficient of a factorizing function
$\wh f^0$ for $f^0$ by (4.17) and $\g$ is a M\"{o}bius automorphism
of $\D$, on which the infima in (4.18) are attained, while by (4.19)
and homogeneity of $J$,
  \be\label{4.22}
r^n = r^n J^0(f^0) = \inf_{\wh f^0}|\wh \beta_n| |\wh c_1^0| r^n +
O(r^{n+1}).
\end{equation}
This implies $\inf_{\wh f^0}|\wh \beta_n| = 1$ and
$$
|\wh c_1^0| = \fc{|c_1^0|}{|(\g^* \kp)^\prime(0)|} = 1.
$$
By Schwarz's lemma, the last equality can hold only when the cover
$\wh f^0(z) = \epsilon z, \ |\epsilon| = 1$, and then by Lemma 3.3
$f^0(z) = \kappa(z)$ up to rotation. This yields also
$$
|c_n^0| = |c_1^0| = |c_1(\g^* \kappa)| = 2/e,
$$
violating Parseval's equality $\sum_0^\iy |c_n|^2 = 1$ for the
boundary function $\kp(e^{i \theta}), \ \theta \in [0, 2 \pi]$ in
(so $|c_n(\kp)| < 2/e$ for all $n > 1)$. This contradiction proves
(4.21).

It follows that the extremals of $J$ must be of the form
 \be\label{4.23}
 f^0(z) = c_0^0 + c_2^0 z^2 + \dots.
\end{equation}
Now, if $n = 2$, comparison of (4.23) with Lemmas 3.3 and 4.4 (for
$m = 2$) and the equalities similar to (3.3), (4.22) implies that
necessarily $|\wh c_2^0| = 1$; hence $\wh f^0(z) = \epsilon z^2$.
Therefore, the maximal value $|c_2^0|$ on $\mathcal B_1^0$ must be
equal to $2/e$ and is attained only on $f^0(z) = \kp(z^2)$ (up to
rotations).

If $n \ge 3$, the same arguments as in the proof of (4.21) based on
the equality (4.16) (for $m = 2$) imply that also the second
coefficient $c_2^0$ of any extremal function $f^0$ for $J(f)$ must
vanish; hence,
  \be\label{4.24}
f^0(z) = c_0^0 + c_3^0 z^3 + c_4^0 z^4 + \dots \ .
\end{equation}
In the case $n = 3$, the relations (4.20), (4.22), (4.24) and Lemmas
3.3 and 4.4 (for $m = 3$) imply, similar to the previous case, that
$\wh f^0(z) = \epsilon z^3$, thus $f^0(z) = k(z^3)$ up to rotations
and $\max_{\mathcal B_1^0} |c_3| = \kp^\prime(0) = 2/e$.

Arguing similarly for $n = 4, 5, \dots$, one derives successively
that for each $n$ the extremal function $f^0$ must be of the form
$f^0(z) = c_0^0 + c_n^0 z^n + c_{n+1}^0 z^{n+1} + \dots$ and
coincide with $\kp(z^n)$ up to rotations, which implies
$\max_{\mathcal B_1^0} |c_n| = 2/e$, completing the proof of the
theorem.

\bigskip
In the case $n = 2$, one can apply the above arguments to more
general functional $J(f) = c_2 + P(c_1)$ given by Proposition 1.2.
Since by (1.2) $c_2(\kp) = 0$, one immediately gets that any
extremal $f^0$ of $J$ satisfies (4.21), hence $|J(f^0)| = |c_2^0| =
2/e$, which implies the estimate (1.4).

\bk
\section{Second proof of Theorem 1.1}

We first establish that $\kp_n(z)$ is the maximizing function in
local setting compatible with Schwarz's lemma, which provides the
assertion of Theorem 1.1.

\bk\noindent $\bf 1^0$. \ Consider more general bounded functionals
on $\mathcal B_1$ of the form
 \be\label{5.1}
J(f) = c_n + F(c_{m_1}, \dots, c_{m_s})
\end{equation}
where $c_j = c_j(f); 1 \le n, m_j$ and $F$ is a holomorphic function
of $s$ variables in an appropriate domain of $\C^s$. We assume that
this domain contains the origin $\mathbf 0$ and that $F, \ \partial
F$ vanish at $\mathbf 0$.

Using the factorization (3.1) and the map $\mathbf k$ generated by
the function (1.2) via Proposition 3.1($b$), we obtain a functional
 \be\label{5.2}
J(\kp \circ \wh f) = \wh J(\wh c_1, \dots, \wh c_n)
\end{equation}
on $\wh f(z) = \wh c_0 + \wh c_1 z + \dots \in H_1^\iy$, and
$\sup_{\mathcal B_1^0} |J(f)| = \sup_{H_1^\iy} |\wh J(\wh f)|$.

Noting that all $\wh f \in H^\iy$ belong to the space $\B$, we
define for $\vp \in A_1(\D), \ \psi \in \B$ the scalar product
 \be\label{5.3}
l_\psi(\vp) = \langle \vp, \psi\rangle = \iint\limits_\D (1 -
|z|^2)^2 \vp(z) \ov{\psi(z)} dx dy,
\end{equation}
As was mentioned above, any linear functional on $A_1(\D)$ is of
such a form. Put
 \be\label{5.4}
\nu_\psi(z) = (1 - |z|^2)^2 \ov{\psi(z)}
\end{equation}
and extend the scalar product (5.3) to all $\mu \in L_\iy(\D)$ and
$\vp \in L_1(\D)$ by $\langle \vp, \mu \rangle = \iint_\D \mu \vp dx
dy$. Then
 \be\label{5.5}
\mu - \nu_\psi \in A_1(\D)^\bot = \{\nu \in L_\iy(\D): \ \langle
\nu, \vp\rangle = 0 \ \text{for all} \ \vp \in A_1(\D)\}
\end{equation}
for any $\mu \in L_\iy(\D)$ extending $l_\psi$. In particular, this
holds for the Hahn-Banach extension $\langle \mu_\psi, \vp\rangle$
of $l_\psi$ having the minimal norm.

We shall need some results from the Teichm\"{u}ller space theory.
Define
 \be\label{5.6}
S_F(\z) = \wh f \Bigl(\fc{1}{\z}\Bigr) \fc{1}{\z^4}.
\end{equation}
These functions are \hol \ on the disk
$$
\D^* = \hC \setminus \{\ov \D\} = \{\z \in \hC = \C \cup \{\iy\}: \
|\z| > 1\}
$$
and belong to the space $\B(\D^*)$ with norm $\|\vp\| = \sup_{\D^*}
(|\z|^2 - 1)^2 |\vp|$. Any $S_F \in \B(\D^*)$ is the Schwarzian
derivative
$$
S_F(\z) = \Bigl(\fc{F^{\prime\prime}(\z)}{F^\prime(\z)}\Bigr)^\prime
- \fc{1}{2} \Bigl(\fc{F^{\prime\prime}(\z)}{F^\prime(\z)}\Bigr)^2
$$
of a locally univalent in $\D^*$ function
$$
w = F(\z) = \z + b_1 \z^{-1} + b_2 \z^{-2} + \dots.
$$
By the Ahlfors-Weill theorem, if $\|S_F\| = 2 k < 2$, then $F$ is
univalent on whole disk $\D^*$ and admits $k$-quasiconformal
extension across the unit circle $\{|z| =1\}$ to $\hC$ with Beltrami
coefficient
$$
\mu_F(\z) = \partial_{\ov{\z}} F/\partial_\z F = - \fc{1}{2} (1 -
|\z|^2)^2 \fc{\z^2}{\ov{\z}^2} S_F\Bigl(\fc{1}{\ov \z}\Bigr).
$$
Every element $\mu \in L_\iy(\D)$ we consider as defined everywhere
on $\C$ with $\mu(\z) = 0$ for $\z \in \D^*$.

The Schwarzians $S_F$ of univalent functions $F$ in $\D^*$ with \qc
\ extensions to $\hC$ form a bounded contractible domain in the
space $\B(\D^*)$ which models the \uTs \ $\T$. Its topologies
generated by the norm of $\B(\D^*)$ and by Teichm\"{u}ller's metric
related to $\|\mu_F\|_\iy$ are equivalent. Hence, for small $r > 0$,
$$
\inf \{\|\mu_\psi\|_\iy: \ \|\psi\|_{\B(\D^*)} = r\} > 0.
$$
The corresponding functions $\psi = \wh f(z) = S_F(1/z^4) z^4$ with
$\|\mu_\psi\|_\iy < r$ form a (connected) domain $U_r$ in $H_1^\iy$
containing the origin.

Now we can formulate our first theorem.

\begin{thm} For any functional $J$ of type (5.1), there exists a
number $r_0(J) > 0$ such that for all $r \le r_0(J)$,
 \be\label{5.7}
\sup_{\mathbf k(U_r)} |J(f)| = \sup_{\mathbf k(U_r)} |c_n| = M_n r,
\quad M_n = \max_{\mathbf B_1} |J(f)|.
\end{equation}
\end{thm}

\bk\noindent{\bf Proof}. As was mentioned above, each function $\psi
\in \B$ determines by (5.3) a linear functional on $A_1(\D)$.
Applying to $\psi$ the reproducing formula
  \be\label{5.8}
 \psi(\zeta) = \fc{3}{\pi} \ \iint\limits_\D  \fc{(1 - |z|^2)^2 \psi(z)}
{(1 - \ov z \zeta)^4} dx dy, \quad \zeta \in \D,
\end{equation}
which is valid for all $\psi$ with $\iint_\D (1 - |z|^2)^2 |\psi(z)|
dx dy < \infty$ (see, e.g. \cite{Ah}), one gets for its derivatives
$$
\psi^{(p)}(0) = \fc{3 \cdot 4 \cdot 5 \dots (p + 4)}{\pi p!} \
\iint\limits_\D (1 - |z|^2)^2 \psi(z) \ov z^p dx dy.
$$
Hence the coefficients $\wh c_p(\wh f)$ of $\wh f = \psi$ are
represented via
 \be\label{5.9}
\wh c_p = \fc{(p + 1)\dots(p+4)}{2 \pi} \ \iint\limits_\D (1 -
|z|^2)^2 \psi(z) \ov z^p dx dy,
\end{equation}
or $\wh c_p = M_p^\prime \ov{\langle z, \psi\rangle}, \ p = 0, 1,
\dots$, with
$$
M_p^\prime = (p + 1)\dots(p+4)/(2 \pi).
$$

Consider the bounded linear transformation
 \be\label{5.10}
\mathcal L: \ \mu \mapsto \psi(\z) = \fc{3}{\pi} \iint\limits_\D \
\fc{\mu(z) dx dy}{(1 - \ov z \z)^4}: \ L_\iy(\D) \to \B(\D).
\end{equation}
It satisfies
  \be\label{5.11}
\langle \vp, \mathcal L \mu\rangle = \iint\limits_\D \vp \ov{\mu} dx
dy, \quad\text{for} \ \ \vp \in A_1(\D),
\end{equation}
and, similar to (5.9) the coefficients of $\psi = \mathcal L \mu$
are given by
 \be\label{5.12}
\wh c_p(\psi) = M_p^\prime \iint\limits_\D \mu(z) \ov z^p dx dy
\end{equation}
which represents the variation of coefficients under varying the
elements $\mu \in L_\iy(\D)$.

Both functionals $J$ and $\wh J$ are extended to all such $\psi$ as
the limits of their values on $f_r(z) = \kp \circ \psi(r z) = c_{0r}
+ c_{1r} z + \dots \in \mathcal B_1^0, \ r \nearrow 1$ and letting $
(c_0, \dots, c_m) = \lim_{r\to 1} (c_{0r}, \dots, c_{m r})$ for
finite collections. Denote these extensions by $J(\mu)$ and $\wh
J(\mu)$.

Our goal is to show that for any extremal function $\psi_0 =
\mathcal L \mu_0$ maximizing $J(\mu)$  on a small ball $U_r =
\{\|\mu\|_\iy < r\}$ (whose existence of $f_0$ follows from
compactness) defines generates a function $f_0 = \kp \circ \wh f_0
\in \mathcal B_1^0$. First of all, we have:

\begin{lem} For small $r > 0$, any extremal $\psi_0 = \mathcal L \mu_0$ on
$U_r$ is orthogonal to all powers $z^p$ with $p \ne 1, 2, \dots, n$,
i.e., for all such $p$,
$$
\langle z^p, \psi_0\rangle  = \langle \mu_0, \ov z^p\rangle = 0,
$$
and therefore, $\psi_0$ is a polynomial
 \be\label{5.13}
\psi_0(z) = t \sum\limits_1^n \wh d_j z^j, \quad |t| = 1.
\end{equation}
\end{lem}

\noindent{\bf Proof}. By a computation, for a fixed $r < 1$ and
$\tau \to 0$, one gets, using the relations (5.11) and (5.12),
 \be\label{5.14}
\begin{aligned}
\max_{\|\mu\|_\iy \le \tau r} \left| \wh J(\wh f) \right| &=
\max_{\|\mu\|_\iy \le \tau r} \left| \sum\limits_{j=1}^n
\fc{\partial \wh J}{\partial \wh c_j} d \wh c_j \right| +
O_n(\tau^2) = \max_{\|\mu\|_\iy \le \tau r} \left| \iint\limits_\D
\ov{\mu(z)} \vp_n(z) dx dy \right|  \\
&= \iint\limits_\D \left| \mu_0(z) \ov{\vp_n(z)} \right| dx dy +
O_n(\tau^2) = \left| \wh J(\tau \mu_0) \right| + O_n(\tau^2),
\end{aligned}
\end{equation}
where $\vp_n$ is a polynomial of the form (5.13) whose coefficients
are determined by the initial coefficients of $\kp$.

Now consider for any fixed $p \ne 1, \dots, n$ the auxiliary
functional
$$
\wh J_p(\mu) = \wh J(\mu) + \xi \wh c_p = \wh J(\mu) + \xi
M_p^\prime \langle \mu, \ov z^p \rangle
$$
with $\xi \in \C$. Then, similar to (5.14),
 \be\label{5.15}
\max_{\|\mu\|_\iy \le r} |\wh J_p(\mu)| = r \iint\limits_\D
|\vp_n(z) + \xi M_p^\prime z^p| dx dy + O_n(r^2), \quad r \to 0,
\end{equation}
and the remainder term estimate is independent of $p$. Using the
known properties of the norm
$$
h_p(\xi) = \iint\limits_\D |\vp_n(z) + \xi z^p| dx dy
$$
following from the Royden\cite{Ro1} and Earle-Kra\cite{EK} lemmas,
one obtains from (5.14), (5.15) that for small $\xi$ there should be
$$
\max_{\|\mu\|_\iy \le r} |\wh J_p(\mu)| = \max_{\|\mu\|_\iy \le r}
|\wh J(\mu)| + r o_p(\xi) + O(r^2 \xi) + O(r^2),
$$
For $\|\mu\|_\iy \le \tau r$ with fixed $r$ and $\tau \to 0$, this
estimate yields
  \be\label{5.16}
\max_{\|\mu\|_\iy \le r} |\wh J_p(\mu)| = \max_{\|\mu\|_\iy \le \tau
r} |\wh J(\mu)| + \tau o_p(\xi) + O(\tau^2 \xi) + O(\tau^2).
\end{equation}
On the other hand, as $\xi \to 0, \ \tau \to 0$,
$$
\begin{aligned}
|\wh J_p(\tau \mu_0)| &= |\wh J(\tau \mu_{\psi_0})| + \Re
\fc{\ov{\wh J(\tau \mu_0}}{|\wh J(\tau \mu_0)|} M_p^\prime \xi
\langle \tau \mu_0, \ov z^p\rangle +
O(\tau^2 \xi^2) \\
&= |\wh J(\tau \mu_0)| + \tau \xi M_p^\prime \langle \tau
\mu_{\psi_0}, \ov z^p\rangle + O(\tau^2 \xi^2)
\end{aligned}
$$
by suitable choices of $\xi \to 0$. Comparison with (5.16) implies
the desired orthogonality $\langle \mu_0, \ov z^p\rangle = 0$.

Substituting into (5.10) the expansion
$$
1/(1 - \ov z \z)^4 = (1 + \ov z \z + \ov z^2 \z^2 + \dots)^4 = 1 + 4
\ov z \z + \dots,
$$
one gets after the term-wise integration the representation (5.13),
completing the proof of the lemma.

We now establish that for small $r$ and some $|t| = 1$, the extremal
$\mu_0$ must be of the form
 \be\label{5.17}
\mu_0(z) = r t \mu_n(z) := r t |\vp(z)|/\vp_n(z),
\end{equation}
where $\vp_n$ is the polynomial of order $n$ given above. We suppose
that (5.17) does not hold and show that this leads to a
contradiction. Without loss of generality, one can take $\mu_0$ to
be an extremal $L_\iy$ function arising in the Hahn-Banach extension
of $l_{\psi_0}$ and satisfying (5.5). Pass to functionals
$$
J^0(f) = J(f)/M_n, \quad \wh J^0(\wh f) = J^0 \circ \mathbf k (f)
$$
mapping $\mathcal B_1^0$ onto the unit disk. The differential of
$\wh J^0$ at $\psi = \wh f = \mathbf 0$ defines a linear operator
$P_n: \ L_\iy(\D) \to L_\iy(\D)$ acting by
$$
P_n(\mu) = \a_n \langle \mu, \vp_n\rangle \mu_0.
$$
Then, in view of our assumption,
$$
\{P_n(t \psi_0): \ |t| < 1\} \subsetneqq \{|t| < 1\};
$$
thus by Schwarz lemma,
 \be\label{5.18}
|P_n(r \mu_0)| = \rho(r) < r.
\end{equation}
Now consider the function
$$
\om_0 = \nu_{\psi_0} - \rho(r) \mu_0
$$
with $\nu_{\psi_0}$ defined by (5.4). We show that $\om_0$
annihilates all functions $\vp \in A_1(\D)$.

Lemma 5.2 and the mutual orthogonality of the powers $z^m, \ m \in
\field Z$, yield that $\langle \om_0, z^p\rangle = 0$ for all $p =
0, 1, \dots$ distinct from $n$. So we have only show that
$$
\langle \om_0, \vp_n \rangle =0, \quad \vp_n = z^n.
$$
Take the conjugate operator
$$
P_n^*(\vp) = \a_n \langle \mu_0, \vp\rangle \vp_n
$$
mapping $L_1(\D^*)$ to $L_1(\D^*)$. It fixes the subspace $\{t
\vp_n: \ t \in \C\}$, and $P_n(\om_0) = \mathbf 0$. Thus from (5.4)
and (5.18), for some $t$,
$$
\langle \om_0, \vp_n\rangle = t \langle \om_0, P_n^* \vp_n\rangle =
t \langle P_n\om_0, \vp_n\rangle = 0,
$$
which means that $\om_0 \in A_1(\D)^\bot$.

It is proven in the theory of extremal \qc \ maps (see, e.g.,
\cite{GL}) that one of the characteristic properties of the extremal
elements $\mu_0 \in L_1(\D)$ (the Beltrami coefficients) for
functionals $l_\psi$ on $A_1(\D)$ arising by their Hahn-Banach
extension to $L_1(\D)$ is
   \be\label{5.19}
\|\mu_0\|_\iy \le \inf \{\|\mu_0 + \om\|_\iy: \ \om \in
A_1(\D)^\bot\}.
\end{equation}
Applying (5.5) and (5.20), one obtains
$$
r \le \|\mu_0 - \om_0\|_\iy = \|\rho_1(r) \mu_0\|_\iy =
\|\rho(r)\|_\iy,
$$
which contradicts (5.18) and proves (5.17).

(In fact, in our case we have more, since the Schwarzians $S_F \in
\B(\D^*)$ generated by functions $f \in B_1^0$ via (5.6) satisfy
$|S_f(\z)| = o(|\z| - 1)^2)$ as $|\z| \to 1$; thus their extremal
$\mu_0$ are of Teichm\"{u}ller form $\mu_0 = k |\vp|/\vp$ with $k =
\const$ and $\vp \in A_1(\D)$; for such $\mu_0$ the inequality
(5.19) is strong).

Finally we have to prove that in (5.13)
$$
\psi_0(z) = r t z^n, \quad |t| = 1.
$$
This is a consequence of the extremality of $\psi_0$ which yields
that the image $\psi_0(\D)$ must be a whole disk $\D_r = \{|w| <
r\}$. To establish such a property of the extremal $\psi_0$, on can
apply the following local existence theorem from \cite{Kr1} which we
present here as

\begin{lem}Let $D$ be a finitely connected domain on the Riemann sphere
$\hC$. Assume that there are a set $E$ of positive two-dimensional
Lebesgue measure and a finite number of points
 $z_1, z_2, ..., z_m$ distinguished in $D$. Let
$\a_1, \a_2, ..., \a_m$ be non-negative integers assigned to $z_1,
z_2, ..., z_m$, respectively, so that $\a_j = 0$ if $z_j \in E$.

Then, for a sufficiently small $\ve_0 > 0$ and $\varepsilon \in (0,
\varepsilon_0)$, and for any given collection of numbers $w_{sj}, s
= 0, 1, ..., \a_j, \ j = 1,2, ..., m$ which satisfy the conditions
$w_{0j} \in D$, \
$$
|w_{0j} - z_j| \le \ve, \ \ |w_{1j} - 1| \le \ve, \ \ |w_{sj}| \le
\ve \ (s = 0, 1, \dots   a_j, \ j = 1, ..., m),
$$
there exists a \qc \ automorphism $h_\ve$ of domain $D$ which is
conformal on $D \setminus E$ and satisfies
$$
h_\ve^{(s)}(z_j) = w_{sj} \quad \text{for all} \ s =0, 1, ..., \a_j,
\ j = 1, ..., m.
$$
Moreover, the \Be \ coefficient $\mu_{h_\ve}(z) = \partial_{\bar z}
h_\ve/\partial_z h_\ve$ of $h$ on $E$ satisfies $\| \mu_{h_\ve}
\|_\iy \le M \ve$. The constants $\ve_0$ and $M$ depend only upon
the sets $D, E$ and the vectors $(z_1, ..., z_m)$ and $(\a_1, ...,
\a_m)$.

If the boundary $\partial D$ is Jordan or is $C^{l + \a}$-smooth,
where $0 < \a < 1$ and $l \ge 1$, we can also take $z_j \in
\partial D$ with $\a_j = 0$ or $\a_j \le l$, respectively.
\end{lem}

Now, assuming that $\psi_0(\D) \ne \D_r$ (and hence $\D_r \setminus
\ov{\psi_0(\D)}$ is an open set), one can pick a disk $E \Subset
\D_r \setminus \ov{\psi_0(\D)}$ and applying Lemma 5.2 vary the
coefficients $\wh c_0, \wh c_1, \dots, \wh c_n$ of $\psi_0$ by \qc \
automorphisms $h_\ve$ of the disk $\D$ conformal on $\D \setminus
\ov E$ so that $h_\ve \circ \psi_0 \in H_1^\iy$ and each coefficient
$\wh c_j(h_\ve \circ \wh f_0)$ ranges over a small neighborhood of
$\wh c_j(\wh f_0)$ for all $j = 0, 1, \dots, n$. By appropriate
choice of $h_\ve$, one gets $|\wh J(h_\ve \circ \wh f_0)| > |\wh
J(f_0)|$ (equivalently, $|J(\kp \circ h_\ve \circ \wh f_0)| >
|J(f_0)|$), contradicting the extremality of $f_0 = \kp \circ
\psi_0$ on the ball $\|\mu\|_\iy < r$ and on its proper convex
subset consisting on $\mu$ with $\psi = \mathcal \mu \in H_1\iy$.

Therefore, the polynomial (5.13) must map the unit circle $\{|z| =
1\}$ onto the circle $\{|w| = r\}$ which is possible only when the
$\psi_0 = r t z^{p_0}$, where $|t| = 1$ and $1 \le p_0 \le n$, in
addition, $p_0$ must divide $n$. Were $p_0 < n = k p_0$, then
$$
f_0(z) = \kappa(tz^{p_0}) = \fc{1}{e} + \fc{2}{e} t^{p_0}z^{p_0} +
\dots + c_n^0 t^k z^n + \dots
$$
with $|c_n^0| \ge 2/e$,  which violates Parseval's equality
$\sum_0^\iy |c_m^0 t^m|^2 = r^2$ for the boundary function $\kp(e^{i
p_0 \theta })$. Hence, $\psi_0(z) = t z^n, \ |t| = 1$.

We have established the existence of $r_0 > 0$ such that, for all $r
\le r_0$, any extremal function of the rescaled functional $J^0(f)$
on $\mathbf k(U_r)$ is of the form $f_0 = (t_1 \kp) \circ (t z^n)$
with $|t| = |t_1| = 1$ (and does not depend on $r$). This sharply
estimates $J(f)$ on $\mathbf k(U_r)$ via (5.7) and completes the
proof of Theorem 5.1.

\bk\noindent $\bf 2^0$. \ Now, to derive the assertion of Theorem
1.1, take $J(f) = c_n(f)$. By Theorem 5.1, we have for $|t| = r \le
r_0$ the sharp bound
 \be\label{5.20}
\sup_{\mathbf k(U_r)} |c_n| = c_n(\kp(t z^n)) = 2 r/e, \quad r =
|t|,
\end{equation}
and the extremal $f_0(z) = t z^n$. This estimates yields, together
with Lemma 3.3 and Schwarz's lemma, that the equality (5.20) must
hold for all $r < 1$ and that the functional $j$ is a defining
function of the disk $\{\kp(t z^n): |t| < 1\}$ as a \Ca \ geodesic
in $\mathcal B_1^0$.

For any other geodesic disk $\{t \wh f: \ |t| < 1\}$ in $H_1^\iy$,
the previous arguments provide the strong inequality
 \be\label{5.21}
|c_n(\kp(t \wh f))| < |c_n(t \kp_n)|.
\end{equation}

In the limit as $r \to 1$, one obtains from (5.20) the desired bound
(1.1), with equality for $f = \kp_n$. It remains to show that no
other extremal functions can appear in the limit case $r = 1$.

Let for some $f \in \mathcal B_1$, we have the equality
  \be\label{5.22}
|c_n(f)| = 2/e = c_n(\kp_n).
\end{equation}
Then all function $f_r(z) = f(r z)$ with $r < 1$ belong to $\mathcal
B_1^0$ and
$$
|c_n(f_r)| = 2 r^n/e = c_n(\kp_{n,r}).
$$
Fix $r < 1$ and represent $f_r$ by (3.1) via $f_r = \kp \circ \wh
f_r$ with corresponding $\wh f_r \in H_1^\iy$. Regarding $\wh f_r$
as a point of its geodesic disk $\{t f_r/\|f_r\|_\iy: \ |t| < 1\}$
corresponding to $t = \|f_r\|_\iy < 1$, one obtains from (5.20) and
(5.21), $|c_n(f_r))| \le |c_n(\kp_n)|t$, with equality only for $f =
\kp_n$. This yields that also (5.22) is valid only for $f = \kp_n$,
completing the proof of Theorem 1.1.

\bk I am thankful to Se\'an Dineen for his remarks.

\medskip
{\small\em{ \leftline{Department of Mathematics, Bar-Ilan
University} \leftline{52900 Ramat-Gan, Israel} \leftline{and
Department of Mathematics, University of Virginia,}
\leftline{Charlottesville, VA 22904-4137, USA}}}

\end{document}